\documentclass[]{elsarticle}

\usepackage{graphicx}

\usepackage{amsmath}
\usepackage{amssymb}
\usepackage{bm}
\usepackage{float}
\usepackage{color}
\usepackage{listings}    
\usepackage{lipsum}
\usepackage{graphicx}
\usepackage{mathtools}
\usepackage{wrapfig}
\usepackage{mathptmx} 
\usepackage{multicol}
\usepackage{tensor}

\biboptions{comma,square}

\journal{journal}

\restylefloat{figure}

\definecolor{dkgreen}{rgb}{0,0.6,0}
\definecolor{MyDarkBlue}{rgb}{0.0,0,0.7}
\definecolor{mygray}{rgb}{0.95,0.95,0.95}
\begin{document}

\begin{frontmatter}

\title{Automatic code generator for higher order integrators}

\author[imf]{Asif Mushtaq}
\ead{Asif.Mushtaq@math.ntnu.no}

\author[ify]{K{\aa}re Olaussen}
\ead{Kare.Olaussen@ntnu.no}

\address[imf]{Institutt for matematiske fag, NTNU, N-7491 TRONDHEIM, Norway}
\address[ify]{Institutt for fysikk, NTNU, N-7491 TRONDHEIM, Norway}

\begin{abstract}
Some explicit algorithms for higher order symplectic integration of a large class of
Hamilton's equations have recently been discussed by Mushtaq \emph{et. al}.
Here we present a Python program for automatic numerical implementation
of these algorithms for a given Hamiltonian, both for double precision and multiprecision computations.
We provide examples of how to use  this program, and illustrate behaviour of both the code generator
and the generated solver module(s).
\end{abstract}


\end{frontmatter}

{\bf PROGRAM SUMMARY}

\begin{small}
\noindent
{\em Manuscript Title: } Automatic code generator for higher order integrators.\\
{\em Authors:} Asif Mushtaq, K{\aa}re Olaussen.\\
{\em Program Title:} HOMsPy: Higher Order (Symplectic) Methods in Python\\
{\em Journal Reference:}                                       \\
{\em Catalogue identifier:}                                    \\
{\em Licensing provisions:} None.                               \\
\emph{Programming language: } Python 2.7.                            \\
\emph{Computer: } PC's or higher performance computers.          \\
\emph{Operating system: } Linux, MacOS, MSWindows.           \\
\emph{RAM: } Kilobytes to a several gigabytes (problem dependent).    \\
\emph{Number of processors used: } 1                            
\emph{Classification: } 4.3 Differential equations, 5 Computer Algebra.\\
\emph{External routines/libraries: }
\texttt{SymPy} library~\cite{SymPy} for generating the code. \texttt{NumPy} library~\cite{NumPy}, and optionally
\texttt{mpmath}~\cite{mpmath} library for running the generated code. The \texttt{matplotlib}~\cite{matplotlib}
library for plotting results.\\
\emph{Nature of problem: }
We have developed algorithms~\cite{AAK1} for numerical solution of Hamilton's equations,
\begin{align}
  \dot{q}^a = \frac{\partial H(\bm{q},\bm{p})}{\partial p_a},\quad
  \dot{p}_a = -\frac{\partial H(\bm{q},\bm{p})}{\partial q^a},
  \quad a=1,\ldots,\mathcal{N}
  \label{HamiltonEquations}
\end{align}
for Hamiltonians of the form
\begin{equation}
 H(\bm{q},\bm{p}) = T(\bm{p}) + V(\bm{q}) = \frac{1}{2} \bm{p}^T M \bm{p} + V(\bm{q}),
 \label{Hamiltonian}
\end{equation}
with $M$ a symmetric positive definite matrix. The algorithms preserve the symplectic
property of the time evolution exactly, and are of orders $\tau^N$ (for $2\le N \le 8$)
in the timestep $\tau$. Although explicit, the algorithms are time-consuming and error-prone
to implement numerically by hand, in particular for larger $N$.\\
\emph{Solution method: }
We use computer algebra to perform all analytic calculations
required for a specific model, and to generate the Python code for
numerical solution of this model, including example programs using that code.\\ 
\emph{Restrictions: }In our implementation the mass matrix is assumed to be equal to the unit matrix,
and $V(\bm{q})$ must be sufficiently differentiable.\\
\emph{Running time: }Subseconds to eons (problem dependent).
See discussion in the main article.
\emph{Program: } Python program can be provided on demand from authors.

\end{small}

\section{Introduction}
\label{intro}
The Hamilton equations of motion \eqref{HamiltonEquations}
play an important role in physics and mathematics. They often
require numerical methods for solution \cite{SanzSerna, McLachlanQuispel, HLWG}.
A well-behaved class of such methods are the \emph{symplectic solvers}, which
preserve symplecticity of the time evolution exactly. One simple
way to construct a symplectic solver is to split the time evolutions
into \emph{kicks},
\begin{align}
   \dot{q}^a = 0, \quad 
   \dot{p}_a = -\frac{\partial V(\bm{q})}{\partial q^a},
   \label{theKicks}
\end{align}
which is straightforward to integrate to give
\begin{align}
 q^a(t+\tau) &= q^a(t),\\ 
 p_a(t+\tau) &= p_a(t) - \tau \frac{\partial V(\bm{q}(t),\bm{p}(t))}{\partial q^a},
\end{align}
followed by \emph{moves},
\begin{align}
   \dot{q}^a = \frac{\partial T(\bm{p})}{\partial p_a} = \sum_b M^{ab} p_b, \quad \dot{p}_a = 0,
   \label{theMoves}
\end{align}
which integrates to
\begin{align*}
 q^a(t+\tau)  &= q^a(t) + \tau \sum_b M^{ab} p_b(t+\tau),\\
 p_a(t+\tau) &= p_a(t) - \tau \frac{\partial V(\bm{q}(t),\bm{p}(t))}{\partial q^a}.
\end{align*}
This scheme was already introduced by Newton~\cite{NewtonPrincipica}
(as more accessible explained by Feynman~\cite{FeynmanPhysicalLaw}).
A symmetric scheme is to make a \emph{kick} of size $\frac{1}{2}\tau$, a \emph{move} of
size $\tau$, and a \emph{kick} of size $\frac{1}{2}\tau$ (and repeating).
This is often referred to as the St{\"o}rmer-Verlet method~\cite{Stormer, Verlet};
it has a local error of order $\tau^3$. The solution provided by
this method can be viewed as the exact solution of a slightly different
Hamiltonian system, with a Hamiltonian $H_{\text{SV}}$ which differ
from \eqref{Hamiltonian} by an amount
proportional to $\tau^2$. For this reason the scheme respects
long-time conservation of energy to order $\tau^2$. It will also exactly preserve
conservation laws due to N{\"o}ther symmetries
which are common to $T(\bm{p})= \frac{1}{2} \bm{p}^T \bm{M} \bm{p}$ and $V(\bm{q})$,
like momentum and angular momentum which are often preserved in
physical models~\cite{Goldstein}.

Recently Mushtaq \emph{et.~al.}~\cite{AAK11, AAK12} proposed
some higher order extensions of the St{\"o}r\-mer-Verlet scheme.
These extensions are also based on the 
\emph{kick-move-kick} idea, only with modified Hamiltonians,
\begin{subequations}
\label{EffectiveHamiltonians}
\begin{align}
   H_1 &\equiv T_{\text{eff}} = \frac{1}{2} \bm{p}^T M \bm{p} + \sum_{k\ge1} T_{2k}(\bm{q},\bm{p}),
   \label{Teff}
   \\
   H_2 &\equiv V_{\text{eff}} = V(\bm{q}) + \sum_{k\ge1} V_{2k}(\bm{q}),\label{Veff}
\end{align}
\end{subequations}
where $T_{2k}$ and $V_{2k}$ are proportional to $\tau^{2k}$.
I.e., the proposal is to replace $V(\bm{q})$ in equation \eqref{theKicks} by $V_{\text{eff}}(\bm{q})$, and
$T(\bm{p})$ in equation \eqref{theMoves} by $T_{\text{eff}}(\bm{q},\bm{p})$. 
The goal is to construct $V_{\text{eff}}$ and  $T_{\text{eff}}$
such that the combined \emph{kick-move-kick} process
corresponds to an evolution by a Hamiltonian $H_{\text{eff}}$
which lies closer to the Hamiltonian $H$ of equation~\eqref{Hamiltonian}.
The difference being of order $\tau^{2N+2}$ when summing terms
to $k=N$ in equations~\eqref{EffectiveHamiltonians}.

One problem with this approach is that $T_{\text{eff}}$ in general
will depend on both $\bm{q}$ and $\bm{p}$; hence the \emph{move}-steps
of equation~\eqref{theMoves} can no longer be integrated explicitly.
To overcome this problem we introduce a generating function~\cite{HLWG}
\begin{equation}
   {G}(\bm{q},\bm{P}; \tau) = 
   \sum_{0\le k\le N} {G}_k(\bm{q},\bm{P})\,\tau^k
   \label{generatingFunctionG}
\end{equation}
such that the transformation $\left( \bm{q}, \bm{p}\right) \rightarrow \left( \bm{Q}, \bm{P}\right)$
defined by
\begin{subequations}
\begin{align}
   p_a &= \frac{\partial  {G}}{\partial q^a},\label{CanonicalTransformation}\\
   Q^a &= \frac{\partial {G}}{\partial P_a},\label{CanonicalTransformation2}
\end{align}
\label{CanonicalTransformations}
\end{subequations}
preserves the symplectic structure exactly, 
and reproduce the time evolution generated by $T_{\text{eff}}$
to order $\tau^N$. Here $Q^a$ is shorthand for $q^a(t+\tau)$, and
$P_a$ shorthand for $p_a(t+\tau)$. Equation \eqref{CanonicalTransformation}
is implicit and in general nonlinear, but the nonlinearity is of
order $\tau^3$ (hence small for practical values of $\tau$).
In the numerical code we solve \eqref{CanonicalTransformation} by
straightforward iteration (typically two to four iterations
in the cases we have investigated).

The rest of this paper is organized as follows:
In section \ref{ExplicitExpression} we introduce compact notation
in which we present the general explicit expressions for
$T_{\text{eff}}(\bm{q}, \bm{p})$, $V_{\text{eff}}(\bm{q})$, and $G(\bm{q},\bm{Q})$.
Because of their compactness
these expressions are straightforward to implement in \texttt{SymPy}.

In section \ref{Examples} we provide examples of how to use
the code generator on specific problems. This process proceeds
through two stages: (i) By providing the potential $V(\bm{q})$
(possibly depending on extra parameters) code for
solving the resulting Hamilton's equations (the solver module)
is generated, and (ii) this solver module is used to analyse the model.
The last stage must of course be implemented by the user,
but an example program which explicitly demonstrates how
the solver module can be used is also
generated during the first stage.

The examples given in subsections \ref{VibratingBeam} (Vibrating beam) and
\ref{parametricAnharmonicOscillator} (One-parameter family of quartic anharmonic
oscillators) have known exact solutions; this makes it easy to check whether the
algorithms behave like expected with respect to accuracy. The example
in subsection \ref{TwoDimensionPendulum} (two-dimensional pendulum) demonstrates
that the program can handle nonpolynomial functions, and that it generates
code which preserves angular momentum exactly.

In section \ref{manyOscillator} we consider a collection of many coupled quartic oscillators.
This is intended as a stress-test of the code, investigating how it behaves
with respect to precision as well as time and memory use for larger and more complex
models. Test of both of the code-generating
process, and of the solver modules generated. We find that the latter continue to
behave as expected with respect to accuracy, but that there is room
for improvement in the area of time and memory efficiency, in particular
for large structured systems.

In section \ref{structure} we present the organization of our code generating
program itself, including diagrammatic representations of its structure. 
We finally summarize our experiences in section \ref{concludingRemarks}.

\section{Explicit expressions} \label{ExplicitExpression}
Compact explicit expressions for the terms of order
$\tau^N$ for $N \in \left\{2, 4, 6\right\}$
in equations \eqref{EffectiveHamiltonians} were given in ~\cite{AAK11, AAK12}.
With the notation
\begin{align}
     &\partial_a \equiv \frac{\partial}{\partial q^a},\quad 
     \partial^a \equiv M^{ab}\partial_b,\quad
     p^a \equiv M^{ab} p_b,\quad D \equiv p_a \partial^a,\quad 
     \bar{D} \equiv (\partial_a V)\partial^a,
     \label{diffAdiffB}
\end{align}
where the \emph{Einstein summation convention}\footnote{An index which occur
twice, once in lower position and once in upper position, are implicitly summed
over all available values. I.e, $M^{ab} \partial_{b} \equiv \sum_{b} M^{ab} \partial_{b}$
(we generally use the matrix $M$ to rise an index from lower to upper position).
} is employed, they are
\begin{subequations}
\begin{align}
    T_2 &= -\frac{1}{12} D^2V \tau^2,\label{T2}\\
    T_4 &=  \frac{1}{720}\left( D^4  - 9 \bar{D} D^2 + 3 D\bar{D} \right)V\tau^4,\\
    T_6 &= -\frac{1}{60480}\big(2\, D^6-40\, \bar{D}D^4 +46\,
      D\bar{D} D^3 -15\,D^2 \bar{D}D^2  \nonumber \\
      & +54\,\bar{D}^2 D^2 - 9\,\bar{D}D\bar{D}D - 42\,D\bar{D}^2 D 
        +12\, D^2\bar{D}^2 \big)V\tau^6\\
    V_2 &= \frac{1}{24} \bar{D}V\tau^2,\\
    V_4 &=\frac{1}{480} \bar{D}^2 V\tau^4,\\
    V_6& = \frac{1}{161280}\left( 17\, \bar{D}^3 - 10\,\bar{D}_3\right)V\tau^6.\label{V6}
\end{align}
\label{TsAndVs}
\end{subequations}
In the last line we have introduced
\begin{equation}
 \bar{D}_3 \equiv (\partial_aV)(\partial_b V)(\partial_c V) \partial^a \partial^b \partial^c.
 \label{diffC}
\end{equation}
By introducing the operator $\mathcal{D} = P^a \partial_a$
the explicit expressions for the coefficients ${G}_k$ can be written
{
\begin{subequations}
\begin{align}
    {G}_0 &= q^a P_a,\\
    {G}_1 &= {\textstyle \frac{1}{2}} P^a P_a,\\
    {G}_2 &= 0,\\
    {G}_3 &= -{\textstyle \frac{1}{12}}{\cal D}^2 V,\\
    {G}_4 &=  -{\textstyle \frac{1}{24}}{\cal D}^3 V,\\
    {G}_5 &= -{\textstyle \frac{1}{240}}\left( 3\,{\cal D}^4+ 3\, \bar{D}{\cal D}^2
      -{\cal D} \bar{D} {\cal D} \right) V,\\
    {G}_6 &= -{\textstyle \frac{1}{720}} \left({2\,\cal D}^5  + 8\,\bar{D}{\cal D}^3
      - 5\, {\cal D}\bar{D}{\cal D}^2 \right) V,\\
    {G}_7 &=-{\textstyle \frac{1}{20160}}\Big({10\,\cal D}^6 + 10\,\bar{D}{\cal D}^4 + 90\,{\cal
        D}\bar{D} {\cal D}^3 - 75\,{\cal D}^2 \bar{D} {\cal D}^2 \nonumber\\
     &\phantom{=-{\textstyle \frac{1}{20160}}}\, +\,18\,\bar{D}^2{\cal D}^2 -3\, \bar{D}{\cal
         D}\bar{D}{\cal D} -14\,{\cal D}\bar{D}^2{\cal D} + 4\,
       {\cal D}^2\bar{D}^2 \Big) V,\\
    {G}_8 &= - {\textstyle \frac{1}{40320}} \Big(3\,{\cal D}^7-87\,\bar{D}{\cal D}^5
      +231\,{\cal D}\bar{D}{\cal D}^4 -133\,{\cal D}^2\bar {D}{\cal
        D}^3+ 63\, \bar{D}^2 {\cal D}^3     \nonumber\\
  &\phantom{= - {\textstyle \frac{1}{40320}}} \,-3\,{\cal D}\bar{D}^2{\cal D}^2
     -21\,{\cal D}^2 \bar{D}^2 {\cal D} +\,4\,{\cal D}^3 \bar{D}^2-63
    \,\bar{D}{\cal D}\bar{D}{\cal D}^2\\
    &\phantom{= - {\textstyle \frac{1}{40320}}}
    +25\,{\cal D}\bar{D}{\cal D}\bar{D}{\cal D}\Big) V.\nonumber
\end{align}
\label{AllGs}
\end{subequations}}
The equations (\ref{TsAndVs},\ref{AllGs}), when used in equations 
(\ref{theKicks},\ref{CanonicalTransformations}), define
the \emph{kick-move-kick} scheme for a general potential $V$.
If one uses all the listed terms the
local error becomes  of order $\tau^9$, and the scheme will
respect long-time conservation of energy to order $\tau^8$.

However, explicit implementation of the numerical
code for a specific potential is rather laborious
and error-prone to do by hand, since the repeated differentiations
(with respect to many variables) and multiplications by lengthy expressions
are usually involved. We have therefore
written a code-generating program in \texttt{SymPy},
wich takes a given potential $V$ as input, perform all the necessary algebra symbolically,
and automatically constructs a \texttt{Python} module for solving one full
\emph{kick-move-kick} timestep. 
It also writes a \texttt{Python} program example
using the module; this example may serve as a template for larger applications.

\section{Examples of code generation} \label{Examples}

The submitted code includes a file \texttt{makeExamples.py}, with various examples
which
demonstrate how the code generator can be used. We go through
these examples in this section; they also provide some illustrations of how the
integrators perform in practical use.

\subsection{Vibrating beam} \label{VibratingBeam}

A simple model for a vibrating beam is defined by the Hamiltonian
\begin{equation}
   H = \frac{1}{2} p^2 - \frac{1}{2} q^2 + \frac{1}{4} q^4.
\end{equation}
A Python code snippet which generates a numerical solver
for this problem is the following:

\lstset{
  language=Python,
  title={\normalsize \textit{Creating a module for solving a vibrating beam}},
  numbers=left,
  morekeywords={True}
}
\begin{lstlisting}
def makeVibratingBeam():
    # Choose names for coordinate and momentum
    q, p   = sympy.symbols(['q', 'p'])
    qvars  = [q]; pvars  = [p]
    # Define the potential in terms of the coordinate
    V = -q*q/2 + q**4/4
    # Create code for a double-precision solver 
    kimoki.makeModules('VibratingBeam', V, qvars, pvars)
\end{lstlisting}
Line 8 generates the files \lstinline!VibratingBeam.py!
and \lstinline!runVibratingBeam.py!. The file \lstinline!VibratingBeam.py!
contains the general double precision solver module for this problem.
A simple use of it in an interactive session is illustrated below:

\lstset{
  language=Python,
  title={\normalsize \textit{Interactive use of the solver module}},
  numbers=left,
  morekeywords={numpy, array, kiMoKi}
}
\begin{lstlisting}
>>> from __future__ import division
>>> from VibratingBeam import *
>>> z = numpy.array([1/2, 5/4])
>>> kiMoKi(z); print z
[ 0.62690658  1.28822851]
>>> kiMoKi(z); print z
[ 0.75756578  1.32399846]
\end{lstlisting}
Here \lstinline!z = [q, p]! is a \texttt{NumPy} array containing the current
state of the solution. Each call of \lstinline!kiMoKi! updates this state
(data from previous timesteps are not kept).

\begin{figure}[h]
\begin{center}
\includegraphics[width=0.71\textwidth]{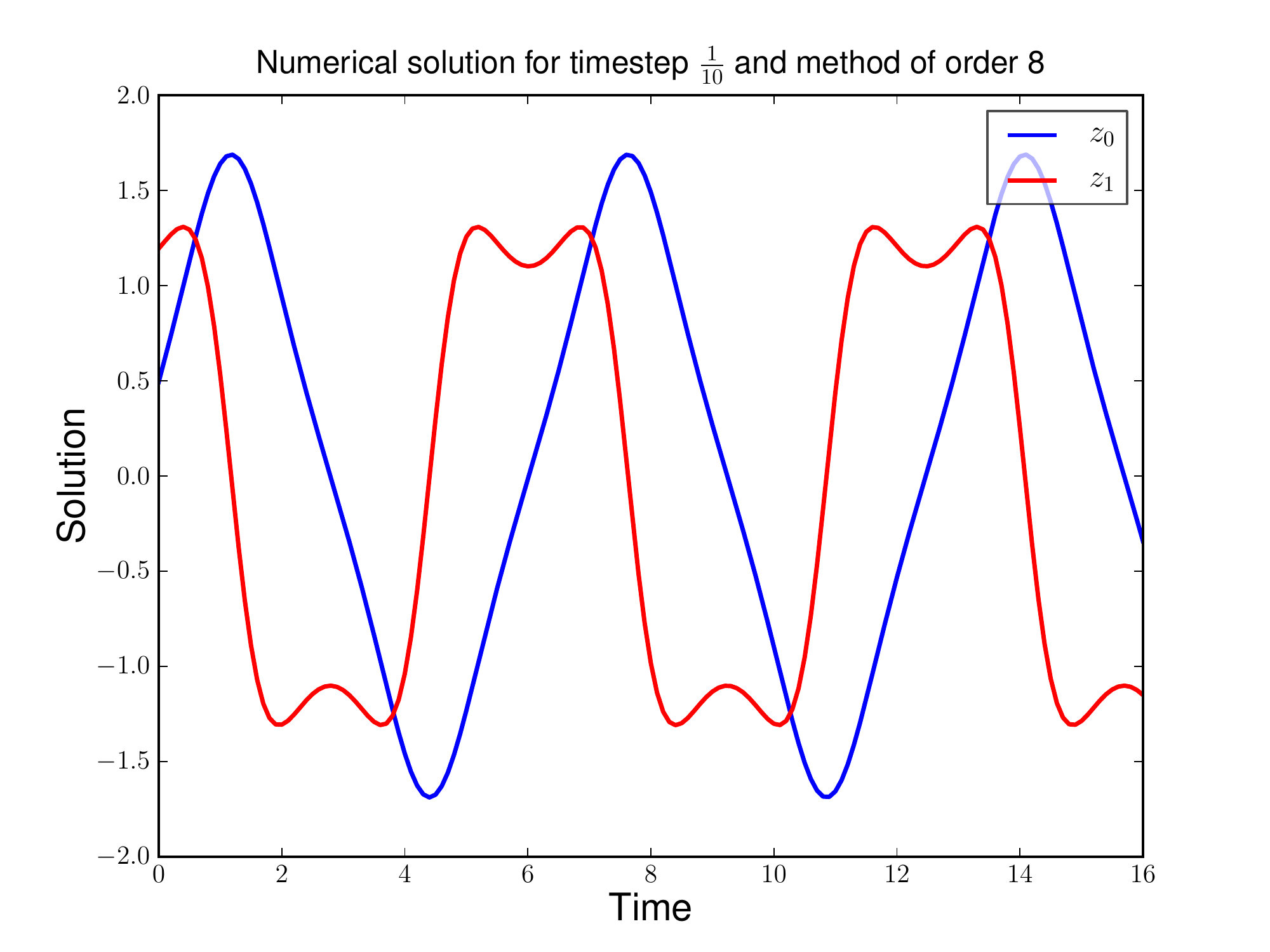}
\end{center}
\caption{Example solution of a vibrating beam with random
initial conditions. Here $z_0$ denotes the position $q$, and $z_1$ the momentum $p$.}
\label{VibratingBeam_SingleSolution}
\end{figure}

The file \lstinline!runVibratingBeam.py! is a small example program demonstrating
basic use of \lstinline!VibratingBeam.py!. A code snippet illustrating
some essential steps is:

\lstset{
  language=Python,
  title={\normalsize \textit{Time evolution of a vibrating beam}},
  numbers=left
}
\begin{lstlisting}
# Select your preferred number of timesteps, order and timestep
    nMax  = 161
    VibratingBeam.order = 8; VibratingBeam.tau   = 1/10
    z = numpy.array([1/2, 5/4]) # Initial condition
    for n in xrange(nMax):
        VibratingBeam.kiMoKi(z) # Integrate one full timestep
\end{lstlisting}
The equation is integrated in line 6. The complete
code in \lstinline!runVibratingBeam.py! is an extension of this snippet.
The initial condition is generated at random, and saved for
possible reuse. Also the complete solution is saved to a (temporary) file
for further processing, together with the parameters \lstinline!tau!,
\lstinline!order!, and \lstinline!nMax!. By running the file
\lstinline!runVibratingBeam.py! a single solution is first generated
and afterwards displayed in a plot. The plot is also saved in
the file \lstinline!VibratingBeam_Soln.png!.
This plot will look similar to Figure \ref{VibratingBeam_SingleSolution}.

\begin{figure}[h]
\includegraphics[width=1.13\textwidth]{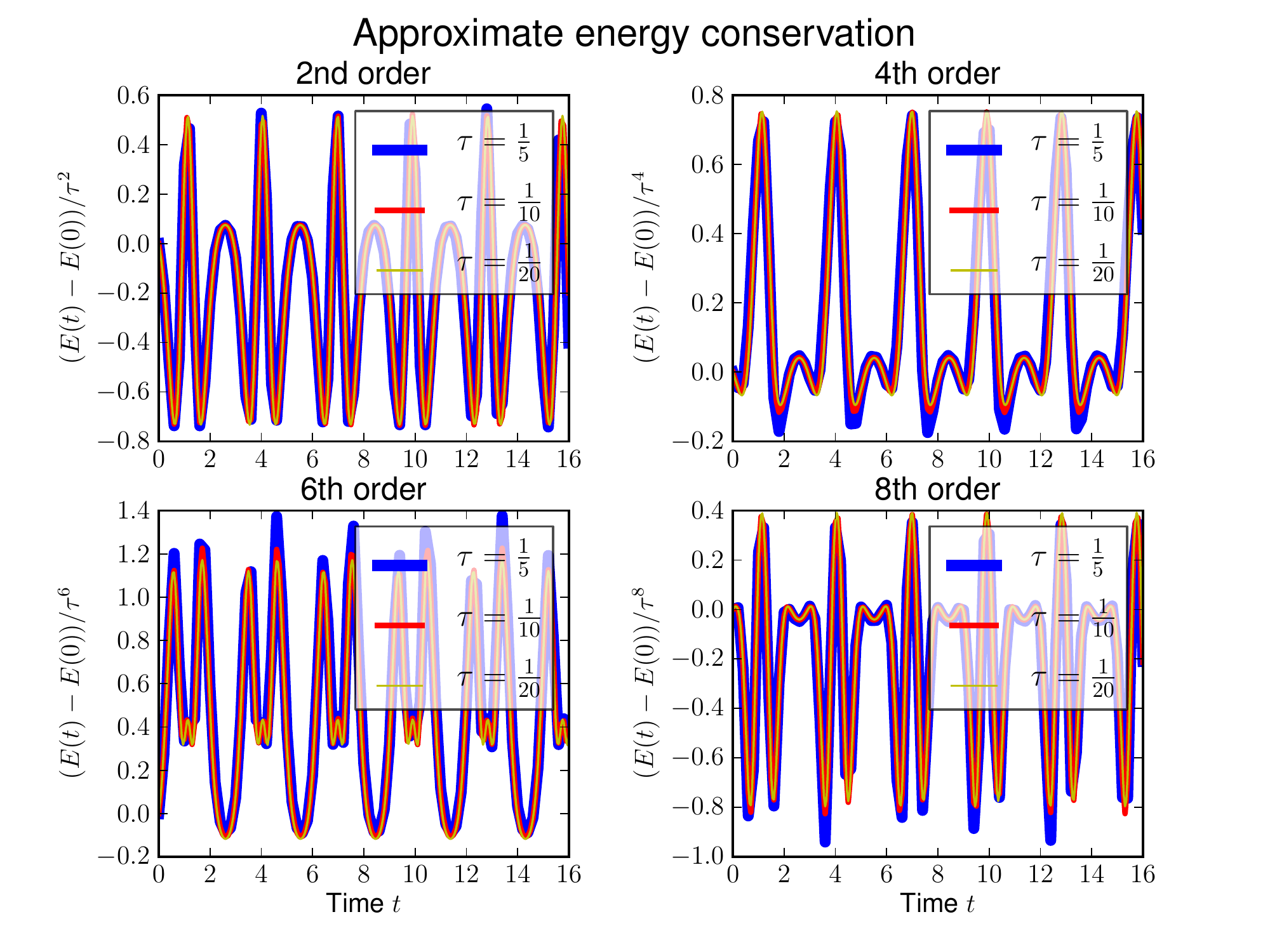}
\caption{Example of how well energy is conserved for integrators of
  various orders $N$,
and how the energy errors scale as expected with the timestep $\tau$. The error decreases
like $\tau^N$ when $\tau$ decreases. In this case, with a periodic solution, there is
also a periodic variation in the energy error.}
\label{VibratingBeam_EnergyErrors}
\end{figure}

To give some impression of the quality of the generated solution,
and how this depends on the timestep and order of the
integrator, the example runfile also make a set of
runs with the same initial condition, but various values of \texttt{tau},
\texttt{order}, and \texttt{nMax}. A simple quality measure,
which is straightforward to implement in general, is
how well the initial energy is preserved as time increases.
This quantity is plotted, with the plot first displayed and
next saved in the file \lstinline!VibratingBeam_EgyErr.png!.
The plot will look similar to Figure \ref{VibratingBeam_EnergyErrors}.
If one prefers to save the plots as \texttt{.pdf}-files the code\\
\lstinline!#import mathplotlib; matplotlib.use('PDF') # Uncomment to ...!\\
on line 17 of the example runfile must be uncommented (then the plot will most likely not be displayed on screen).
During the run process solution data is saved to several \texttt{.pkl}-files; these
are normally removed after the data has been plotted. To keep this data the code\\
\lstinline!os.remove(filename) # Comment out to keep datafile!\\
on lines 108 and/or 202 of the example runfile must be commented out.

As illustrated by Figure \ref{VibratingBeam_EnergyErrors}, and verified by all
other cases we have investigated, the energy error scales like $\tau^N$,
where $\tau$ is the timestep and $N$ is the \texttt{order} of the integrator.
The energy error does not grow with time, but varies in a periodic manner --- 
following the periodicity of the generated solution.
Note that the integration module, in this case \lstinline!VibratingBeam!, contains a
parameter \lstinline!epsilon! which governs how accurate \eqref{CanonicalTransformation}
is solved. We have observed a systematic growth in the energy error
when this parameter is chosen too large, thereby violating
symplecticity (too much). 

Another property of interest and importance is how the average time per integration
step varies with the \lstinline!order! of the method. The run example prints a measure
of the computer time used. The results of this, for a longer run than the unmodified
run example, are shown in
Table \ref{VibratingBeam_IntegrationTime}. As can be seen, the penalty of using
a higher order method is quite modest for a simple model like the vibrating beam,
in particular when the timestep $\tau$ is small.

\begin{table}
\begin{center}
\begin{tabular}{|c|c|c|c|}
\hline
&&&\\[-2ex]
$N$&$\tau=0.2$&$\tau=0.1$&$\tau=0.05$\\[0.2ex]
\hline
&&&\\[-2ex]
2&0.26&0.26&0.26\\
4&1.03&0.84&0.73\\
6&1.23&1.01&0.88\\
8&1.54&1.25&1.09\\[0.2ex]
\hline
\end{tabular}
\end{center}
\caption{Evaluation time per integration step in milliseconds,
dependent on timestep $\tau$ and method \texttt{order} $N$.
The higher order methods ($N>2$) run faster for smaller $\tau$
because the iterative solution of \eqref{CanonicalTransformation}
requires fewer iterations, hence becomes faster. This example, like all others,
have been run on a workstation equipped with two four-core Intel Xeon E5520
processors. The \emph{ratio} between numbers like those above are probably more
relevant that their absolute values.
}
\label{VibratingBeam_IntegrationTime}
\end{table}

\begin{figure}[h]
\includegraphics[width=1.0\textwidth]{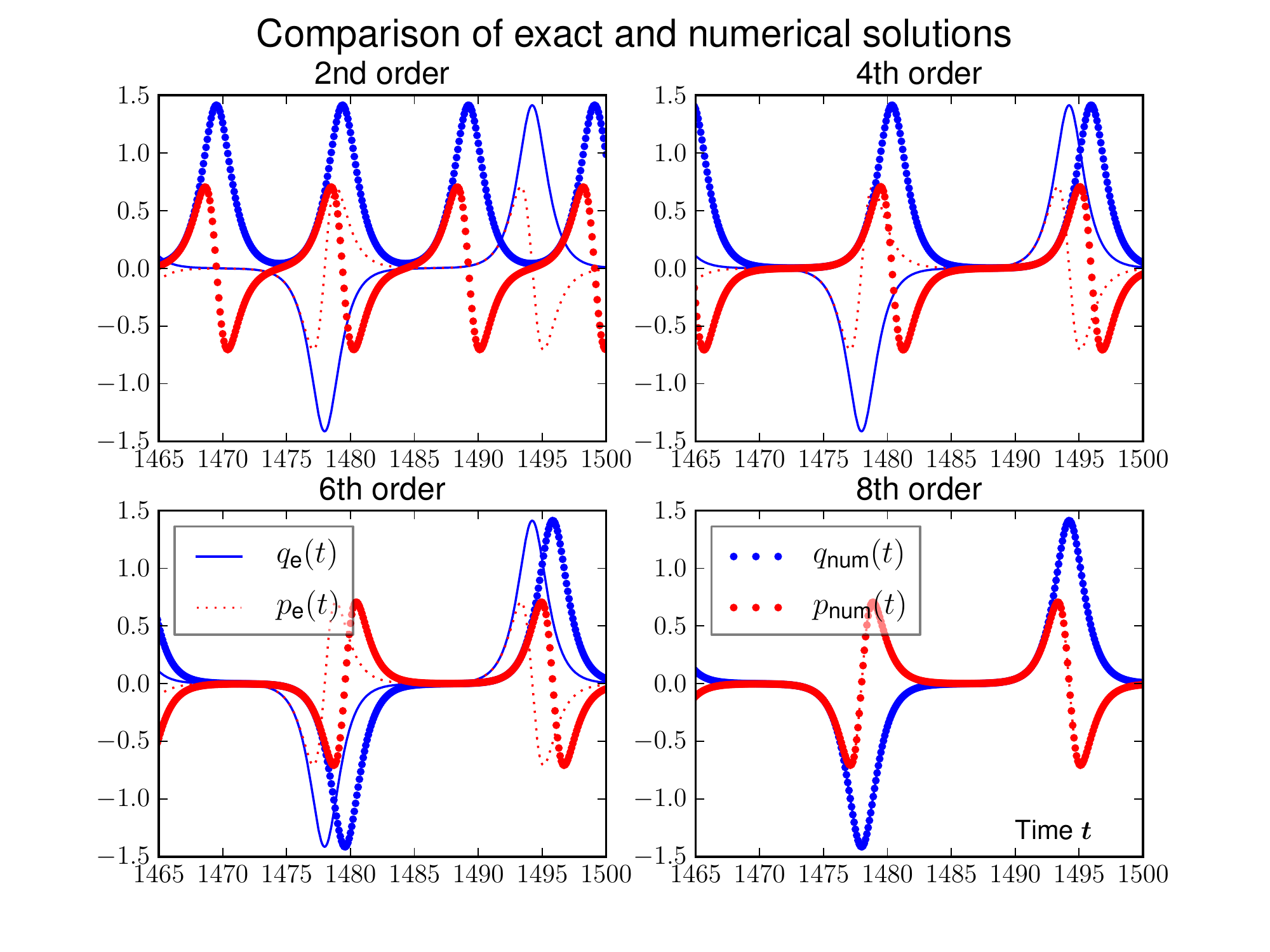}
\caption{Comparison of the exact elliptic solution \eqref{ExactEllipticSolution}
with those found by our automatically generated numerical code, for the initial condition
$q_0=\sqrt{2} + 10^{-6}$ and $p_0 = 0$, with timestep $\tau=10^{-1}$.
The thick blue [red] lines show the numerical solutions
for $q(t)$ [$p(t)$], the thin blue [red] lines show the exact
solution for  $q(t)$ [$p(t)$].
}
\label{ComparisonOfAnharmonicOscillator}
\end{figure}

The exact solution of this problem can be expressed in terms of
Jacobi elliptic functions, cf.~equations~\eqref{ExactEllipticSolution} below.
This allows direct comparison between the exact and the numerical solutions,
as shown in Figure \ref{ComparisonOfAnharmonicOscillator}.
The initial condition is chosen such that the energy is close to the critical energy,
$E=0$, where the solutions bifurcates from motion over the potential
hill at $q=0$ to motion in only one of the
two potential wells. The exact solution moves over the potential hill.
As can be seen, this is respected by the solutions of order $N=6$ and $8$,
but not by the solutions of order $N=2$ (St{\"o}rmer-Verlet) and $4$. This
demonstrates that there may be cases where a higher order method, or an impractically
small stepsize $\tau$, is required to obtain even the qualitatively correct solution.

Up to the time $t$ we have computed and plotted the solutions,
the order $N=8$ numerical solution cannot be visually distinguished
from the exact one in this plot. A difference would become visible for sufficiently
large $t$, because the two solutions have slighly different periods.
The difference in periods is of order $\tau^8$ for $N=8$.

\subsection{One-parameter family of quartic anharmonic oscillators}  \label{parametricAnharmonicOscillator}

A generalization of the previous problem is the class of non-linear oscillators defined by the Hamiltonian
\begin{equation}
\label{AnharmonicOscillators}
   H = \frac{1}{2} p^2 + \frac{1}{2} \alpha \,q^2 + \frac{1}{4}q^4.
\end{equation}
A one-parameter class of exact solutions to this problem can be expressed in terms of
Jacobi elliptic functions~\cite{AbramowitzStegun},
\begin{subequations}
\label{ExactEllipticSolution}
\begin{align}
   q(t) &=   q_0\,\text{cn}(\nu t \vert \text{k}),\\
   p(t) &=  -q_0 \nu\,\text{sn}(\nu t \vert \text{k})\,\text{dn}(\nu t \vert \text{k}).
\end{align}
\label{AllSolutions}
\end{subequations}
Here the initial conditions are
$q(0) = q_0$, and $p(0) =0$. The vibrating beam discussed
in the previous subsection corresponds to the case of $\alpha=-1$. 
Equations \eqref{AllSolutions}
exhaust the set of solutions which have a maximum
at $q_0 >0$ (or a minimum at $q_0 < 0$). This condition imposes the
restriction that $\alpha + q^2_0 \ge 0$. The parameters and
energy of the solution are
\begin{equation}
   \nu = \left(\alpha + q^2_0\right)^{1/2},\qquad
   k   = 2^{-1/2}\,q_0/\nu,\qquad
   E   = \frac{1}{2} \alpha\, q_0^2 + \frac{1}{4} q_0^4.
\end{equation}
For $\alpha + \frac{1}{2}q^2_0 > 0$ the energy is positive, and $q(t)$
oscillates symmetrically around $q=0$; for $\alpha + \frac{1}{2}q^2_0 < 0$
the energy is negative, and $q(t)$ oscillates in one of the two possible
potential wells (depending on the sign of $q_0$).
For $\alpha + q^2_0 < 0$ the solution has a minimum at $q_0 >0$ 
(or maximum at $q_0 < 0$).\footnote{In which case the solution can be written
\begin{equation}
   q(t) = q_0\,\text{nd}(\nu t, k),\qquad
   p(t) = \nu q_0\,k^2\,\text{sd}(\nu t, k)\,\text{cd}(\nu t, k),
\end{equation}
with $\nu = \left(-\alpha - q^2_0/2\right)^{1/2}$ and 
$k = \nu^{-1}\left(-\alpha - q_0^2\right)^{1/2}$. 
In \eqref{AllSolutions}
the modulus $k>1$ when $\alpha + q_0^2/2 <0$. An alternative expression,
with $0 < k < 1$, is
\begin{equation}
   q(t) = q_0\,\text{dn}(\nu t, k),\qquad
   p(t) = -\nu q_0\, \text{sn}(\nu t, k)\,\text{cn}(\nu t, k),
\end{equation}
where $\nu= 2^{-1/2} q_0$ and  $k= q_0^{-1}\left[2(\alpha + q^2_0)\right]^{1/2}$.
}

A code snippet for generating numerical solvers for this problem is the following
\lstset{
  language=Python,
  title={\normalsize \textit{Solving a one-parameter class of anharmonic oscillators}},
  numbers=left
}
\begin{lstlisting}
def makeAnharmonicOscillator():
    # Choose names for coordinate, momentum and parameter
    q, p, alpha  = sympy.symbols(['q', 'p','alpha'])
    qvars = [q]; pvars = [p]; params = [alpha]
    # Define potential in terms of coordinate and parameters
    V = alpha*q**2/2 + q**4/4
    # Code for double-precision and multiprecision computations 
    kimoki.makeModules('AnharmonicOscillator', V, qvars, pvars, 
                        PARAMS=params, MP=True, VERBOSE=True)    
\end{lstlisting}

The code in line 8-9 shows that the \lstinline!makeModules! function
may take optional arguments: If the Hamiltonian depends on a
list of parameters, this list must be assigned to the keyword
\lstinline!PARAMS!. If the \lstinline!MP! keyword is set
to \lstinline!True! then two additional files are generated:
In this case the files \texttt{AnharmonicOscillatorMP.py},
which is a solver module using multiprecision
arithemetic, and \texttt{runAnharmonicOscillatorMP.py},
which is a runfile example using this multiprecision solver.
When the \lstinline!VERBOSE! keyword is set to \lstinline!True!
some information from the code generating process will
be written to screen, mainly information about the time used
to process the various stages. This may be of use when the
code generation process takes a very long time, as will happen
with complicated models.

\begin{figure}[H]

\includegraphics[width=0.90\textwidth]{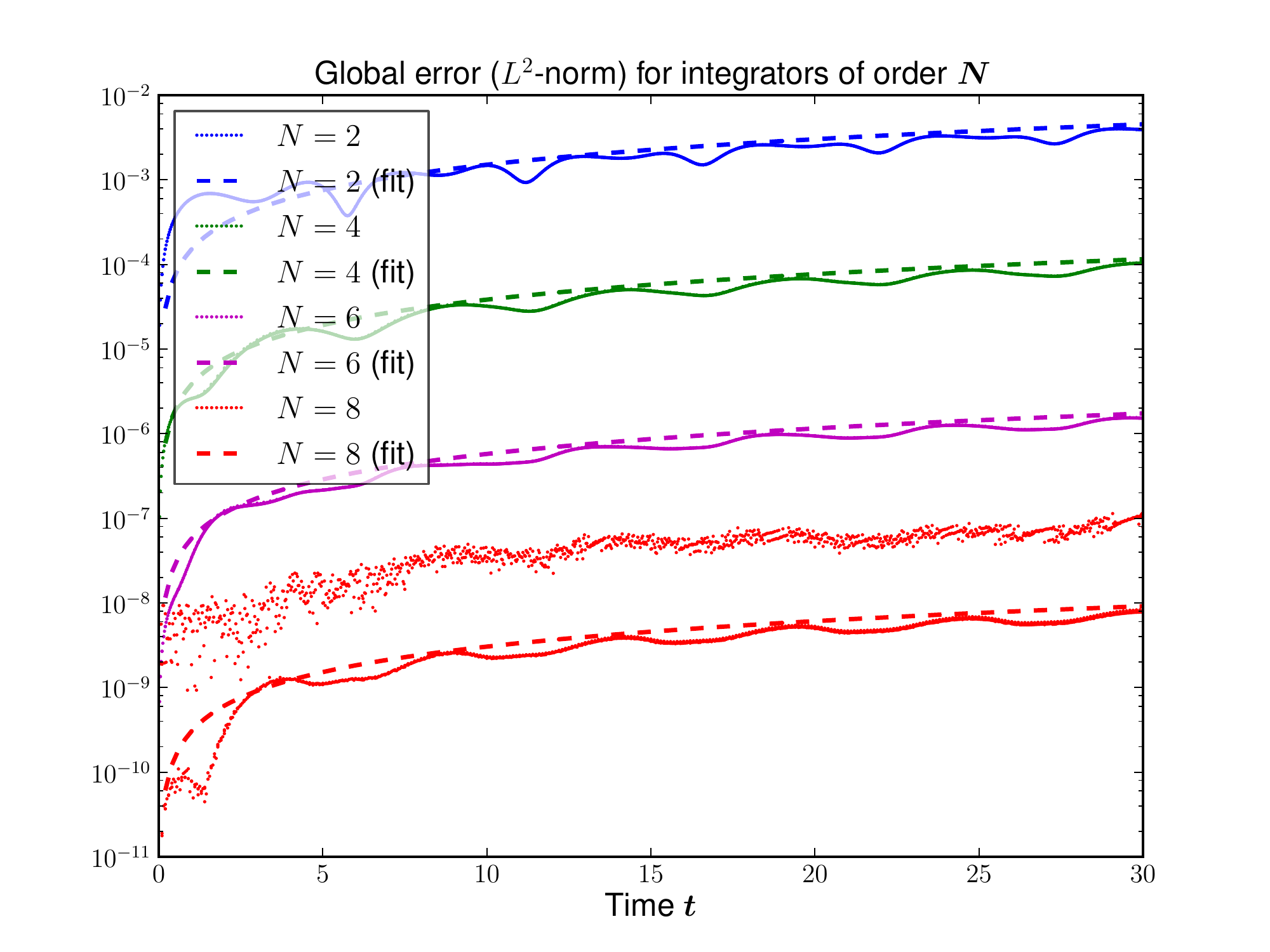}
\caption{The global error \eqref{GlobalError} computed for integrators of various
orders $N$, and stepsizes $\tau = \frac{1}{5}$, $\frac{1}{10}$,
$\frac{1}{20}$, $\frac{1}{40}$. To check the $\tau^N$-scaling the
errors for the last three $\tau$'s are multiplied by respectively
$2^{N}$, $4^{N}$, $8^{N}$. The dashed lines are crude fits to
equation \eqref{GlobalErrorFit} with $C_N$ as the fitting parameters.
The plot is not untypical, with the $N=8$ result deviating
from \eqref{GlobalErrorFit} for the smallest values of $\tau$.
This is a consequence of the finite accuracy of double precision
calculations, not a failure of the $N=8$ integrator.
}
\label{GlobalL2Error}
\end{figure}

\begin{figure}[H]

\includegraphics[width=0.90\textwidth]{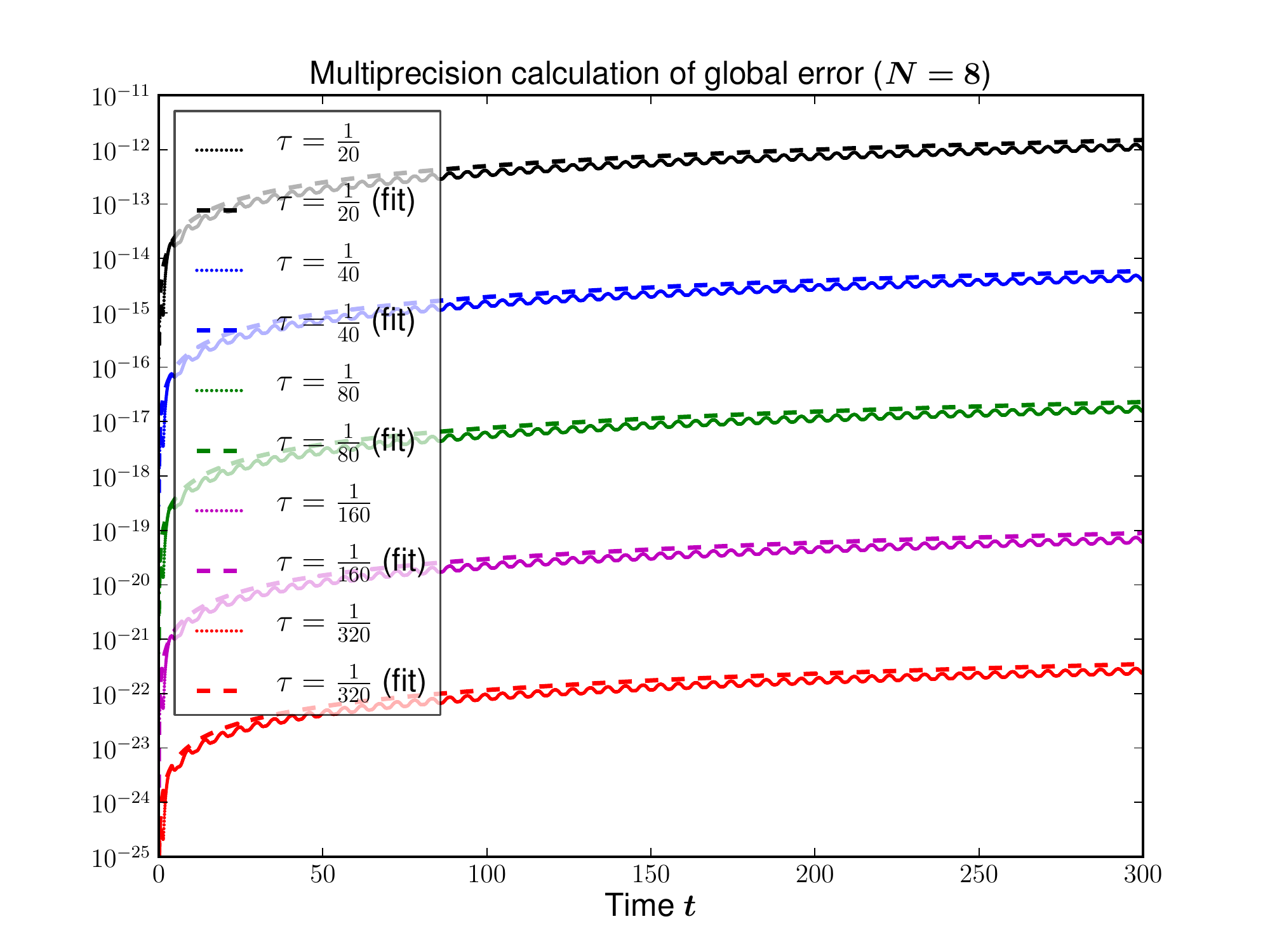}
\caption{The global error \eqref{GlobalError} computed to 35 decimals precision
for the $N=8$ integrator, with stepsizes $\tau = \frac{1}{20}$, $\frac{1}{40}$,
$\frac{1}{80}$, $\frac{1}{160}$, and $\frac{1}{320}$.
The dashed lines correspond to equation \eqref{GlobalErrorFit},
with $C_8$ fitted crudely to the $\tau=\frac{1}{20}$ data.
}
\label{GlobalL2ErrorMP}
\end{figure}

To check the accuracy of the numerical solution in more detail,
we have modified \texttt{runAnharmonicOscillator.py} to
\texttt{analyseAnharmonicOscillator.py}, where the \emph{global error}
\begin{equation}
    \epsilon(t) \equiv \| (q_{\text{(e)}}(t_m),p^{\text{(e)}}(t_m))
    - (q_{\text{(n)}}^m, p^{\text{(n)}}_m) \|
    \label{GlobalError}
\end{equation}
is computed (for random values of $\alpha$ and $q_0$), and plotted.
Here the sub$|$super-script ${\text{(e)}}$ labels the exact solution, and ${(\text{n})}$ the
numerical one. One resulting plot, for parameters $\alpha=0.13$
and $q_0=0.54$, is shown in Figure \ref{GlobalL2Error}. In general,
the global error fits well to the formula, cf.~theorem 3.1 in the book \cite{HLWG},
\begin{equation}
    \epsilon(t)\equiv C_N \,\tau^N\;t,
    \label{GlobalErrorFit}
\end{equation}
where $N$ is the order of the integrator, and $C_N$ is independent
of $\tau$ but depends on the parameters of the model, the initial conditions,
and the norm $\|\cdot\|$ used in \eqref{GlobalError}. We have used the $L^2$
norm in Figures \ref{GlobalL2Error} and \ref{GlobalL2ErrorMP}. 

It can be deduced from Figure \ref{GlobalL2Error} that
to fully exploit the power of the higher order integrators
one must go beyond double precision accuracy. We have therefore
implemented an option (\lstinline!MP=True!) for automatically generating
multiprecision versions of the integrators and
run examples. The file \texttt{analyseAnharmonicOscillatorMP.py}
is an adaption of \texttt{runAnharmonicOscillatorMP.py} which
computes the global error to multiprecision accuracy. The result
for $N=8$ and various small values of $\tau$ (and the same parameters
as in Figure \ref{GlobalL2Error}) is shown in Figure \ref{GlobalL2ErrorMP}.

\subsection{Two-dimensional pendulum} \label{TwoDimensionPendulum}

As a final example of this section we want to demonstrate that our program
can handle non-polynomial potentials as well.
Hence we consider the Hamiltonian for a slightly
distorted\footnote{The motion of a real pendulum is constrained
to the surface of a sphere, which cannot be described
by a constant mass matrix.}
version of a two-dimensional pendulum
\begin{equation}
   H = \frac{1}{2}\left(p_0^2 + p_1^2\right) -\cos\left(\sqrt{q_0^2 + q_1^2}\right).
\end{equation}
Here both the kinetic and potential energy is invariant under rotations;
hence we expect the generated code to preserve angular momentum,
\begin{equation}
    L(t) = \left[q_0(t) p_1(t) - q_1(t) p_0(t)\right],
\end{equation}
exactly.

\lstset{
  language=Python,
  title={\normalsize \textit{Solving a two-dimensional pendulum}},
  numbers=left
}
\begin{lstlisting}
def makeTwoDPendulum():
    # Choose names for coordinate, momentum and parameter
    q0, q1, p0, p1 = sympy.symbols(['q0', 'q1', 'p0', 'p1'])
    qvars = [q0, q1]; pvars = [p0, p1]
    # Define potential in terms of coordinate and parameters
    V = -cos(sqrt(q0**2+q1**2))
    # Code for multiprecision computation only 
    kimoki.makeModules('TwoDPendulum', V, qvars, pvars, DP=False,
                       MP=True, MAXORDER=6, VERBOSE=True)    
\end{lstlisting}
Here we demonstrated one additional optional argument of \lstinline!makeModules!, 
\lstinline!MAXORDER!, which  can be used the restrict the maximum order 
of solvers being generated (6 in this example). 

\begin{figure}[H]
\begin{center}
\includegraphics[width=0.99\textwidth]{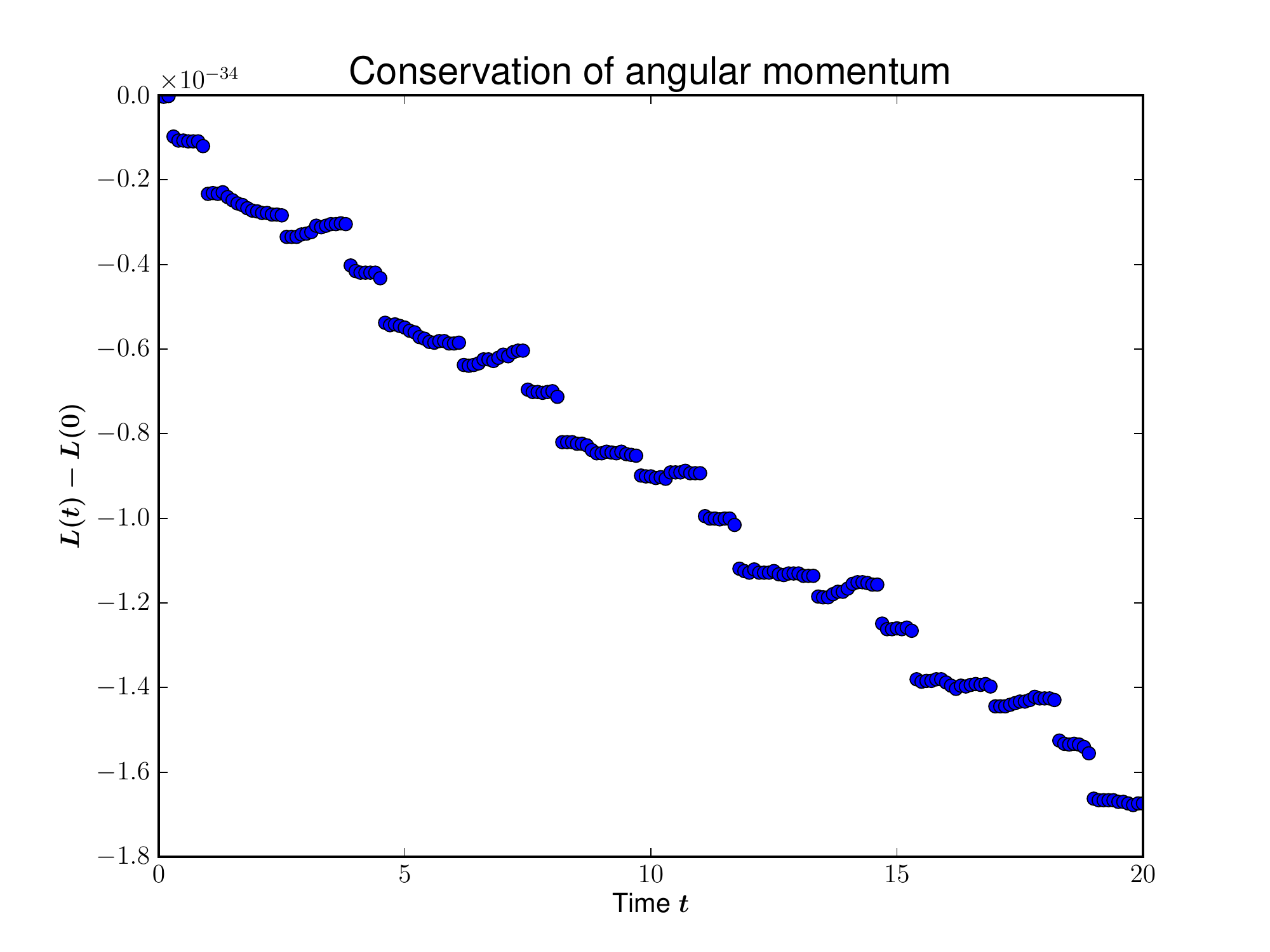}
\caption{This figure demonstrates that our integration algorithm respects
conservation of a angular momentum exactly (within numerical accuracy).
Here we have integrated Hamilton's equations with a timestep \lstinline!tau = 1/mpf(10)!
with an \lstinline!order = 6! multiprecision solver, using 35 decimal
digits accuracy. Obviously, this figure only displays how roundoff errors
are accumulated with time.
}
\label{TwoDPendulum_AngularMomentum.pdf}
\end{center}
\end{figure}

\section{Analysis of many anharmonic oscillators} \label{manyOscillator}

Consider a sum of Hamiltonians like \eqref{AnharmonicOscillators},
\begin{equation}
   H = \sum^{\cal N}_{a=1} \left[\frac{1}{2} P^2_a + 
     \frac{1}{2}\alpha_a\,\left(Q^{a}\right)^2 + 
     \frac{1}{4}\,\left(Q^{a}\right)^4\right].
   \label{SimpleHamiltonian}
\end{equation}
Since the corresponding Hamiltonian equations of motion decouple,
the solution for each pair $(Q^a, P_a)$ is given by expressions
like \eqref{ExactEllipticSolution}. A direct numerical solution of
this model would not provide any additional test of the integrators.
However, if we make an orthogonal coordinate transformation,
\begin{equation}
    Q^a = \sum^{\cal N}_{j=1} \tensor{R}{^a_j}\,q^j,\qquad
    P_a = \sum^{\cal N}_{j=1} \tensor{R}{_a^j}\,p_j,
\end{equation}
we obtain an expression
\begin{equation}
    H = \frac{1}{2}\sum^{\cal N}_{j=1} p^2_j + V(\bm{q}),
    \label{ComplicatedHamiltonian}
\end{equation}
where $V(\bm{q})$ \emph{looks} like a general polynomial
potential in $\cal N$ variables with quadratic
and quartic terms.
We expect the numerical algoritms to behave like the
general case for this model, while the exact solution is known in the form
\begin{equation}
    q_{\text{(e)}}^j = \sum_a \tensor{R}{_a^j}\,Q_{\text{(e)}}^a,\qquad
    p^{\text{(e)}}_j = \sum_a \tensor{R}{^a_j}\,P^{\text{(e)}}_a.
\end{equation}

\begin{figure}[H]
\begin{center}
\includegraphics[clip, trim=10.5ex 7ex 12ex 6ex, width=0.81\textwidth]{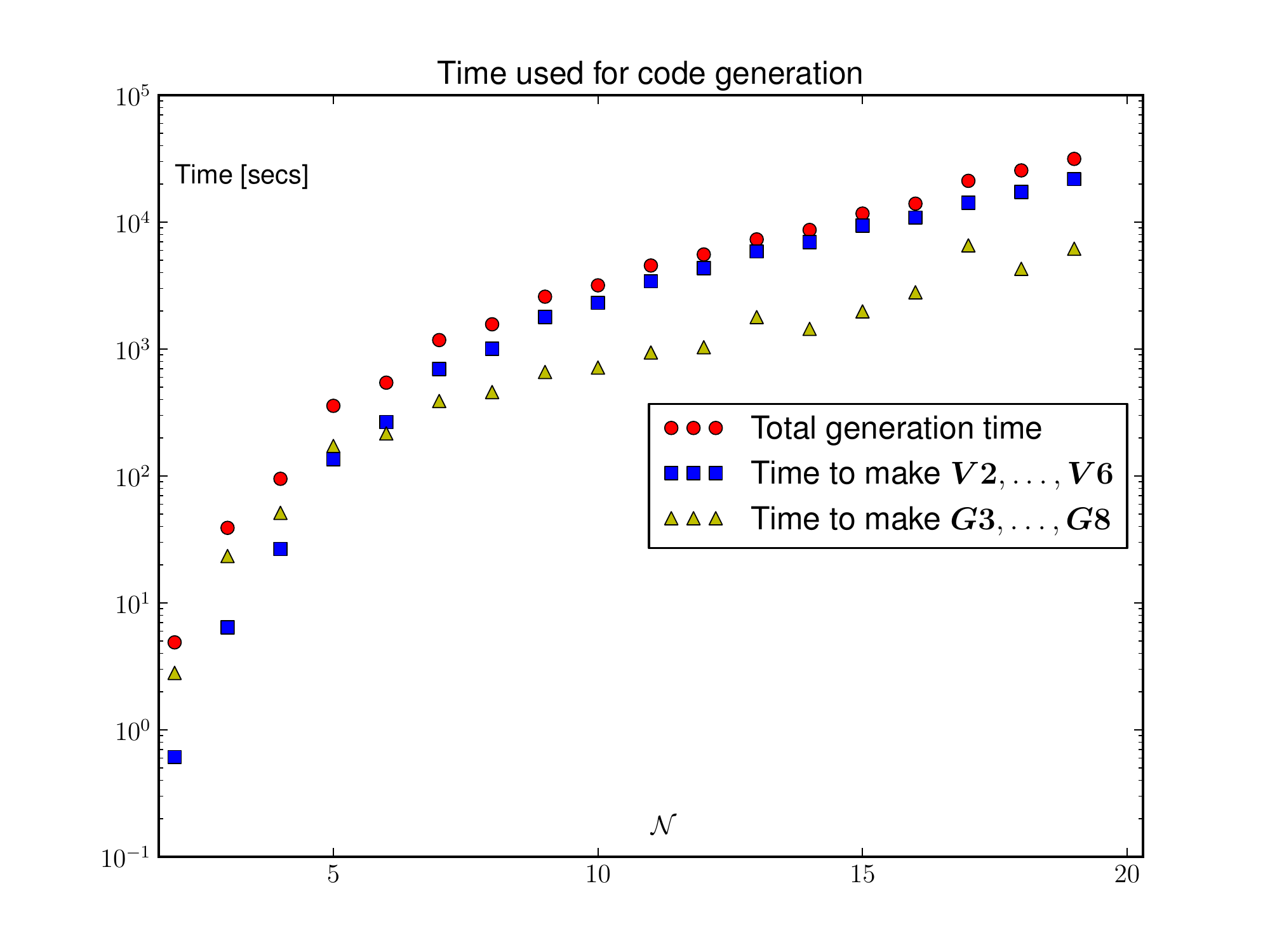}
\caption{This figure illustrates how the time used for code generation scale
with system size $\cal{N}$ for the Hamiltonian \eqref{ComplicatedHamiltonian}.
When ${\cal N}$ is large the most time-consuming individual task seems to be the
computation of $\bar{D}_3 V$.
}
\label{generationTime}
\end{center}
\end{figure}

We have generated such Hamiltonians, using a random orthogonal matrix
with rational coefficents $\tensor{R}{^a_j}=\tensor{R}{_a^j}$, for a range of ${\cal N}$-values
(the way we construct $\tensor{R}{^a_j}$  only nearest- and next-nearest-neighbor
couplings are generated between the variables $q_j$).
This allows us to investigate how the code generator behave for models of
increasingly size and complexity. As illustrated in
Figure \ref{generationTime} the time used to generated the solver module
increases quite rapidly with the number ${\cal N}$ of variables. 

\begin{figure}[H]
\begin{center}
\includegraphics[clip, trim=5.3ex 7ex 11ex 5.5ex, width=0.49\textwidth]{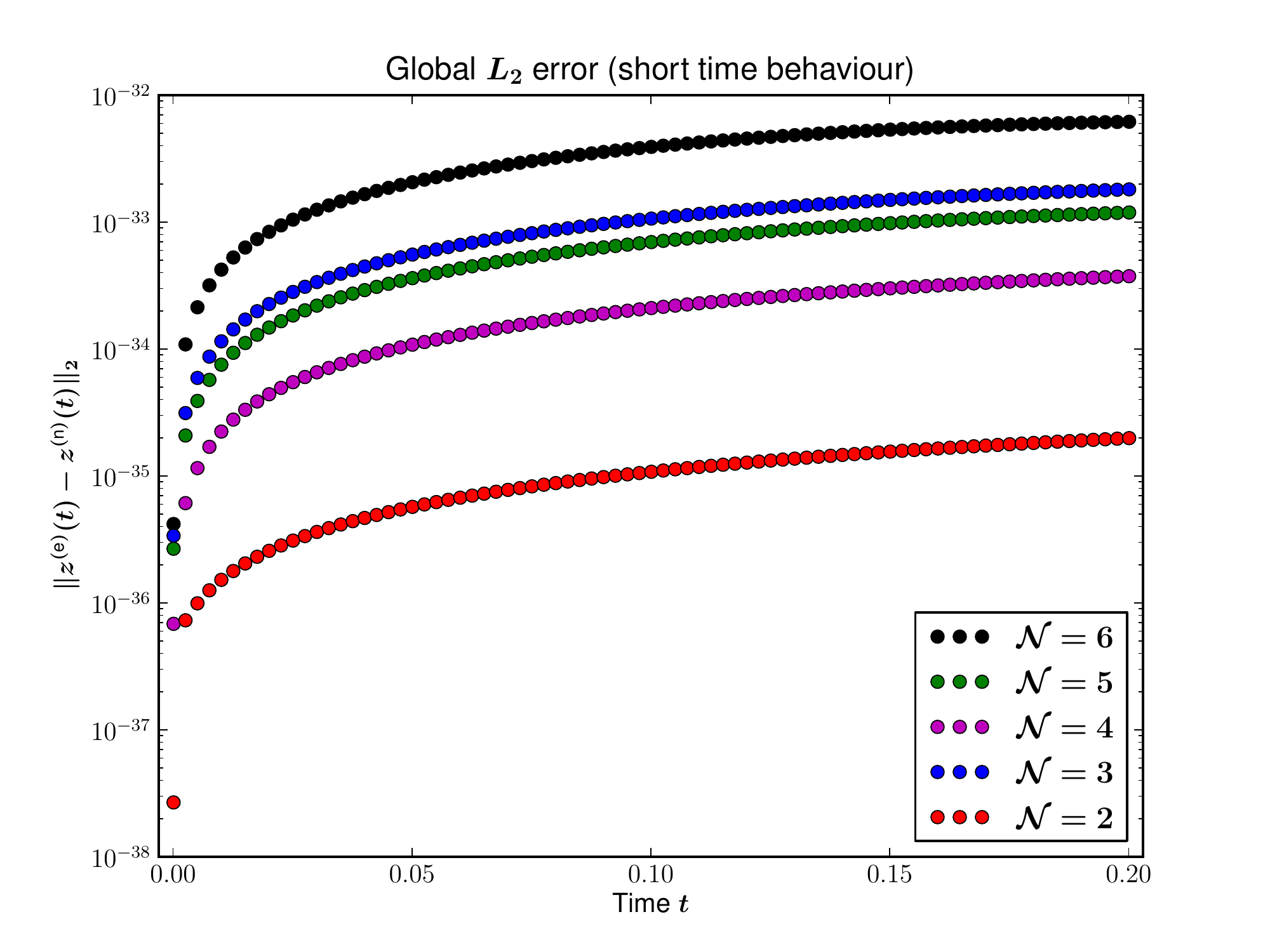}
\includegraphics[clip, trim=5.3ex 7ex 11ex 5.5ex ,width=0.49\textwidth]{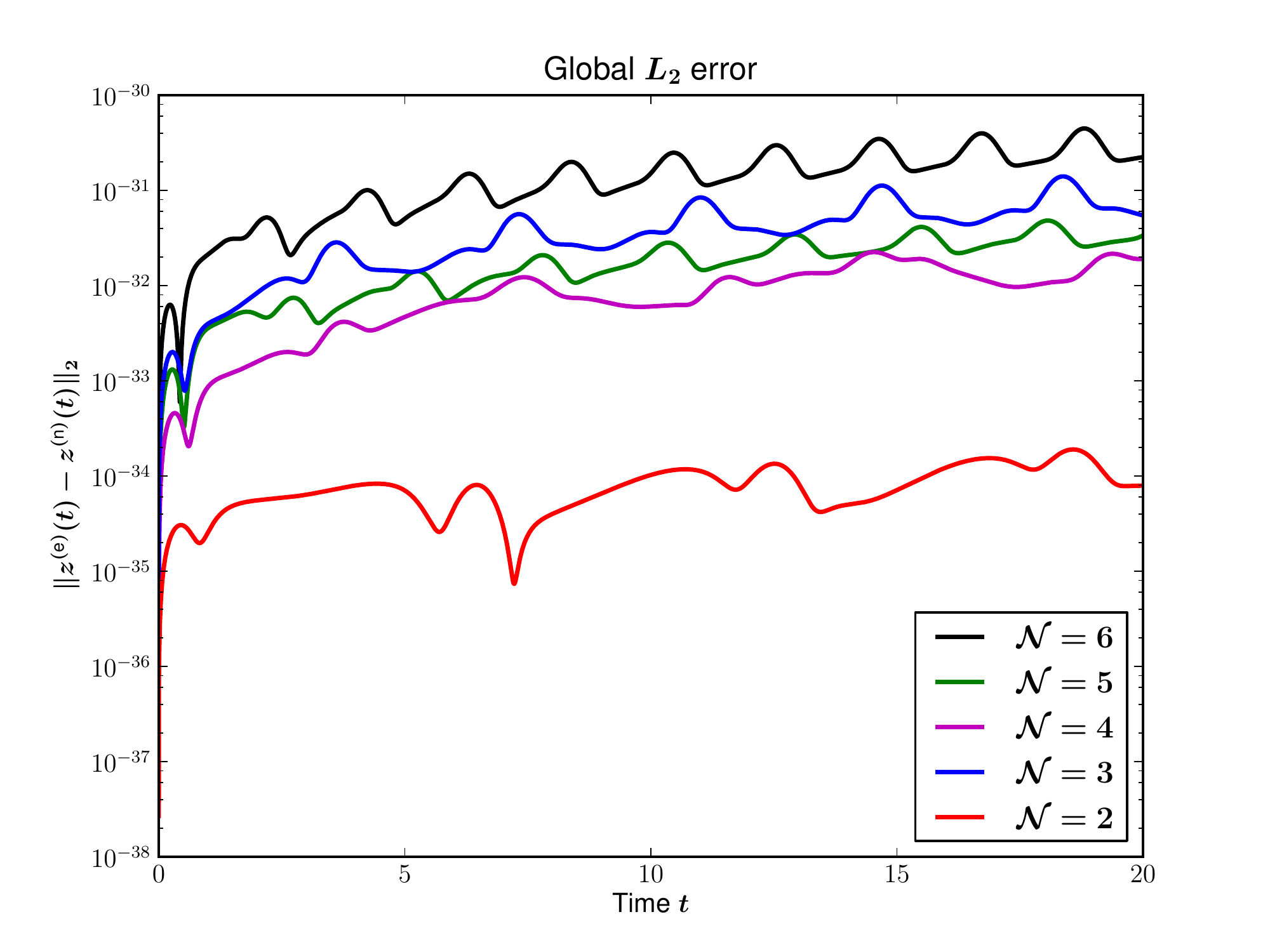}
\caption{
These figures illustrates how the global error varies with time, system size $\mathcal{N}$,
and initial conditions. The short-time behavior is shown to the left; the long-time
behaviour, for the same model and initial conditions, to the right.
We have used a multiprecision integrator of order 8,
with $\tau=10^{-3}$ and computations to 35 decimals precision. 
}
\label{ManyOscillatorGlobalL2Error}
\end{center}
\end{figure}

In Figure \ref{ManyOscillatorGlobalL2Error} we illustrate how the global error
in these models behave. Although there is a general trend that the accuracy
detoriates with system size, this trend is not strictly followed (as can
be seen by the case of $\mathcal{N}=3$). This is a reflection of the fact
that both the models and their initial conditions are generated with a
certain degree of randomness.


\section{Structure of the programs} \label{structure}

\begin{center}
\begin{figure}[H]
\begin{center}
\includegraphics[width=0.8\textwidth]{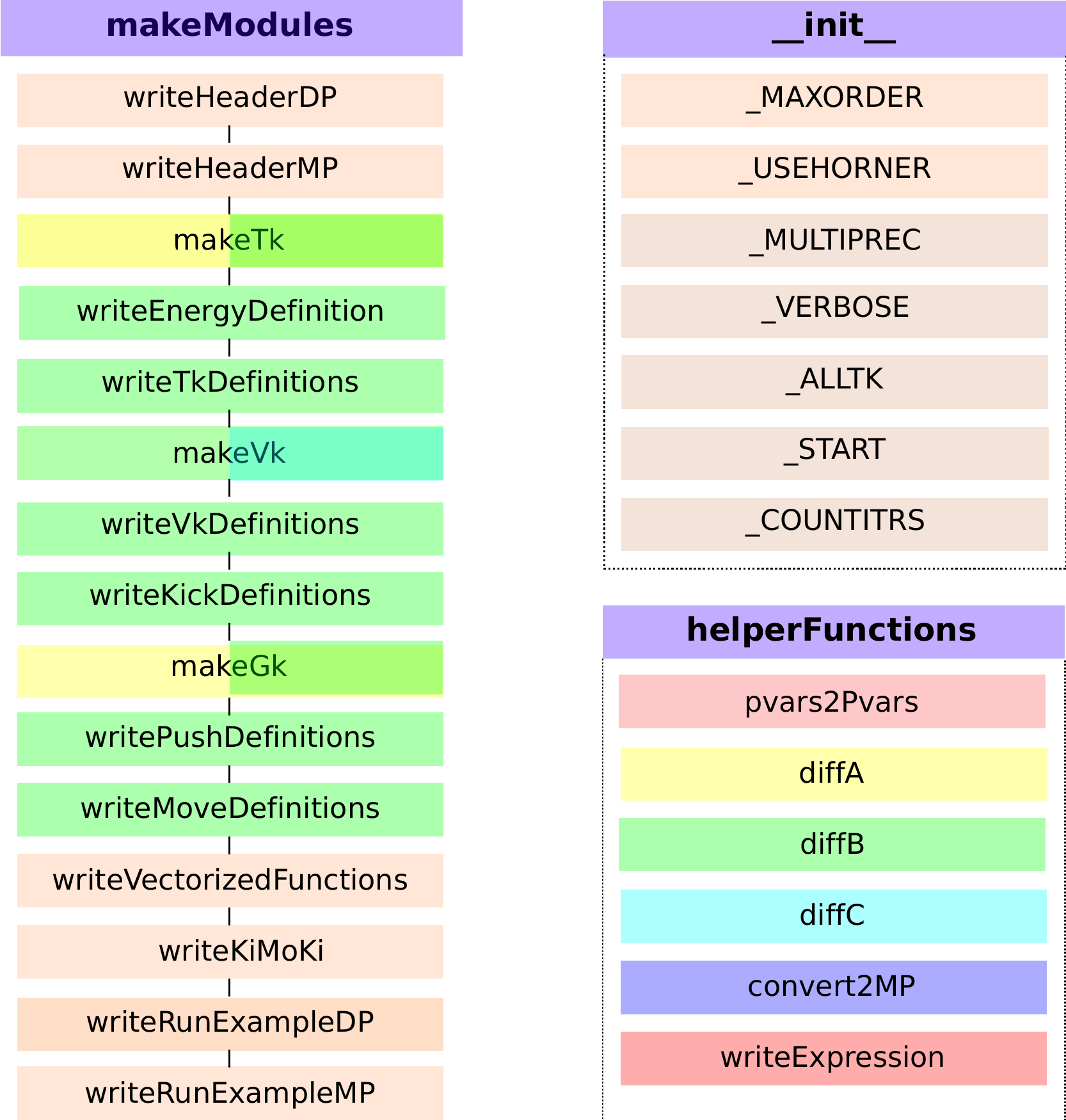}
\end{center}
\caption{Components of the code generating program. The main process is run
in the \lstinline!makeModules(...)! routine;
this calls a sequence of subroutines as illustrated in the left column.
This routine takes a number of optional keyword arguments,
the values of some of them are stored in global variables (as listed in
the \lstinline!__init__! box). Some functions, most of them required at several
places in the program, are collected in \lstinline!helperFunctions!. 
All routines listed in the left column, plus \lstinline!__init__! and
\lstinline!helperFunctions! are defined in files with the same name plus the extension
\lstinline!.py!. These files are in turn collected in the directory (folder)
\lstinline!kimoki!, the name of the code-generating module.
}\label{programFlow}
\end{figure}
\end{center}

From a functional point of view our code generating program
can be characterized as a \emph{module}, hence it should (to our understanding)
be organized into a single file. However, with $1700+$ lines of code and comments
this would make code development and maintainance impractical.
Hence we have organized it as a \emph{package}. I.e., as a set of
\texttt{.py}-files, including a file named \texttt{\_\_init\_\_.py},
in a directory (folder) with the same name as the module (\lstinline!kimoki!). A graphical overview
of the main components of this structure is illustrated
in Figure~\ref{programFlow}.

For a given Hamiltonian this module generates
one or two \emph{solver modules} with routines for numerical solution
of the model, each with a \emph{runfile example} using the solver.
Each runfile example is intended to demonstrate and check basic properties
of its solver module, and to be modified into more useful programs by the user.

In a separate \texttt{examples} directory there is
a file named \texttt{makeExamples.py}, containing the examples
discussed in section \ref{Examples}. By running \texttt{makeExamples.py}
the \texttt{examples} directory will, after some time,
be populated with several solver modules and runfile examples.
Running the runfile examples will in turn generate many \texttt{.pkl}-files with
numerical data (which are normally deleted after use), and some \texttt{.png}-files with plots
of the solution, and how well the solution respects energy conservation. 

\subsection{Diagramatic overview of the code generating module}

\subsection{Brief description of the routines}

\begin{itemize}

\item[$\circ$] \lstinline!makeModules(modname, V, qvars, pvars[, kwargs])!

This is the main subroutine, and the only one intended to by called
by the user.
Here \lstinline!modname! is a ``basename'' of the generated files,
\lstinline!V! is a symbolic expression for the potential $V$,
\lstinline!qvars! is a list of symbolic positions variables (generalized coordinates),
and \lstinline!pvars! is a list of symbolic momentum variables (canonically conjugate momenta).
This routine takes a number of optional keyword arguments (\lstinline!kwargs!) with defaults:

\begin{itemize}
  \item[] \lstinline!PARAMS=None! is a list of symbolic parameters used in $V$.
  \item[] \lstinline!MAXORDER=8! is the maximum order of the generated solvers.
  \item[] \lstinline!DP=True! is a switch which determines if code for
    normal (double) precision solvers should be generated.
  \item[] \lstinline!MP=False! is a switch which determines if code for multiprecision solvers should be generated.
  \item[] \lstinline!VERBOSE=False! is a switch which determines if messages from 
the code generation process should be written.
  \item[] \lstinline!ALLTK=False! is a switch wich determines if code for calculating
    $T_2$, $T_4$, and $T_6$ is generated.
  \item[] \lstinline!COUNTITRS=False! is a switch which determines if code for monitoring the iterative solution of
    \eqref{CanonicalTransformation} is generated. This code makes a histogram of the number of iterations used for solutions.
\end{itemize}

\item[$\circ$] \lstinline!diffA(V, qvars, pvars)!, \lstinline!diffB(V, V0, qvars)!, \lstinline!diffC(V, qvars)!

These three functions implement respectively the operator $D$ defined in equation \eqref{diffAdiffB},
the operator $\bar{D}$ defined in expressions~\eqref{diffAdiffB}, and the operator
$\bar{D}_3$ defined in equation \eqref{diffC}. Here \lstinline!V0! is a symbolic expression
of the potential defining the model, while \lstinline!V! can be any symbolic expression depending
on \lstinline!qvars!, \lstinline!pvars!, and potential parameters \lstinline!params!. The
function \lstinline!diffA! is used by the routines \lstinline!makeTk!
and \lstinline!makeGk!; the function \lstinline!diffB! is used by the
routines \lstinline!makeTk!, \lstinline!makeVk! and \lstinline!makeGk!; the function
\lstinline!diffC! is used by the routine \lstinline!makeVk!.

\item[$\circ$] \lstinline!writeExpression(outfile, expression)!, \lstinline!convert2MP(match)!

The function \lstinline!writeExpression! writes \lstinline!expression! to
\lstinline!outfile! in python format. If a multiprecision version of the solver
module is generated, the function \lstinline!convert2MP! is used to assure that
fractions are converted to multiprecision format. Optionally a polynomial
\lstinline!expression! can be converted to Horner form first (this is not
recommended if \lstinline!expression! depends on many variables, due to time
and memory use).

\item[$\circ$] \lstinline!writeHeaderDP(outfile, pvars, params)!,\\
\lstinline!writeHeaderMP(outfile, pvars, params)!

These routines write the header part of respectively the
double precision and multiprecision solver modules. Some important
solver module variables are defined here, with defaults:

\begin{itemize}

\item[] \lstinline!tau = 1/10! (DP), \lstinline!tau = 1/mpf(1000)! (MP).

The timestep $\tau$ used by the solvers.

\item[] \lstinline!epsilon = 1/10**12! (DP), \lstinline!epsilon = 1/mpf(10**20)! (MP).

The accuracy to which \eqref{CanonicalTransformation} must be solved.
We have observed that the symplectic preserving property of the solvers
is lost when \lstinline!epsilon! is too large, but it must be somewhat
larger than the numerical precision used.

\item[] \lstinline!order = MAXORDER!. Which order of solver to use, setting \lstinline!order! larger
than \lstinline!maxorder! (see below) has no effect.

\item[] \lstinline!params!.  A list of the symbolic potential parameters; these parameters
must be set before starting a solution. 

\item[] \lstinline!maxorder = MAXORDER!. The maximum order of generated solvers. 
Must not be changed by the user.

\item[] \lstinline!dim!. The number of phase space variables. Must not be changed.

\item[] \lstinline!itrs[20]!. A histogram of how many iterations are used to solve \eqref{CanonicalTransformation}.
Exists only if \lstinline!COUNTITRS! is set to \lstinline!True!. Must not be changed by the user.

\end{itemize}

\item[$\circ$] \lstinline!makeTk(V0, qvars, pvars)!

Calculates the contributions $T_0$, $T_2$, $T_4$, and $T_6$ to $T_{\text{eff}}$, cf.~\eqref{Teff},
using explicit expressions in \eqref{TsAndVs}. $T_0$ is used
by the routine \lstinline!writeEnergyDefinition!; $T_2$, $T_4$,
and $T_6$ by the routine \lstinline!writeTkDefinitions!. $T_2$, $T_4$, and $T_6$ are computed
only when the optiononal parameter \lstinline!ALLTK=True!, otherwise they are set to 0.

\item[$\circ$] \lstinline!writeEnergyDefinition(outfile, qvars, pvars, T0, V0)!

Writes the definition of the function \lstinline!energy(z)!, which evaluates
the energy $T_0(\bm{p}) + V_0(\bm{q})$, cf.~equation \eqref{Hamiltonian},
at the phase space point $z \equiv (\bm{q}, \bm{p})$.

\item[$\circ$] \lstinline!writeTkDefinitions(outfile, qvars, pvars, Tk)! 

Writes the definitions of the functions \lstinline!T2(z)!, \lstinline!T4(z)! and \lstinline!T6(z)!,
using the symbolic expressions in the list \lstinline!Tk! calculated by \lstinline!makeTk!.

\item[$\circ$] \lstinline!makeVk(V0, qvars)!

Calculates the contributions $V_2$, $V_4$, and $V_6$ to $V_{\text{eff}}$, cf.~equation
\eqref{Veff}, using explicit expressions in equation \eqref{TsAndVs}. $V_2$, $V_4$, and
$V_6$ are used by the routines \lstinline!writeVkDefinitions!
and \lstinline!writeKickDefinitions!.

\item[$\circ$] \lstinline!writeVkDefinitions(outfile, qvars, Vk)!

Writes the definitions of the functions \lstinline!V2(q)!, \lstinline!V4(q)! and \lstinline!V6(q)!,
using the symbolic expressions in the list \lstinline!Vk! calculated by \lstinline!makeVk!.

\item[$\circ$] \lstinline!writeKickDefinitions(outfile, qvars, Vk)!

Calculates the symbolic expressions $-{\partial V_{\text{eff}}}/{\partial q^a}$,
cf.~equation \eqref{theKicks}, using symbolic expressions in the list \lstinline!Vk!
calculated by \lstinline!makeVk!. These expressions are used to define the
functions \lstinline!kick!$a$\lstinline!(z)! used in the \emph{kick}-steps
of the solvers.

\item[$\circ$] \lstinline!makeGk(V0, qvars, Pvars)!

Calculates the contributions $G_3$, $G_4$, \ldots, $G_8$ to the generation function 
$G(\bm{q}, \bm{P}; \tau)$, cf.~equation \eqref{generatingFunctionG}, using the explicit
expressions in equation \eqref{AllGs}. $G_3$, $G_4$, \ldots, $G_8$ are used by
\lstinline!writePushDefinitions! and \lstinline!writeMoveDefinitions!.

\item[$\circ$] \lstinline!writePushDefinitions(outfile, qvars, Pvars, Gk)!

Calculates the symbolic expressions ${\partial G}/{\partial q^a}$,
cf.~equation \eqref{CanonicalTransformation}, using symbolic expressions in the list \lstinline!Gk!
calculated by \lstinline!makeGk!. These expressions are used to define the
functions  \lstinline!push!$a$\lstinline!(z)! used in the \emph{push}-steps
of the solvers.

\item[$\circ$] \lstinline!writeMoveDefinitions(outfile, qvars, Pvars, Gk)!

Calculates the symbolic expressions ${\partial G}/{\partial P_a}$,
cf.~equation \eqref{CanonicalTransformation2}, using symbolic expressions in the list \lstinline!Gk!
calculated by \lstinline!makeGk!. These expressions are used to define the
functions  \lstinline!move!$a$\lstinline!(z)! used in the \emph{move}-steps
of the solver(s).

\item[$\circ$] \lstinline!writeVectorizedDefinitions(outfile, qvars, Pvars)! 

Writes definitions of functions \lstinline!vecKicks(idx, z)!, \lstinline!vecPushes(idx, z)!,
and \lstinline!vecMoves(idx, z)!, \emph{simulating} parallel evaluation of
\lstinline!kick!$a$\lstinline!(z)!, \lstinline!push!$a$\lstinline!(z)!,
and \lstinline!move!$a$\lstinline!(z)!, for all $a$ in the list \lstinline!idx!.
We don't think parallel evaluation is actually achieved\footnote{This depends on
how the \lstinline!numpy.vectorize! function is implemented.}; hence presently this
is only \emph{syntactic sugar} simplifying the main solver routine, \lstinline!kiMoKi(z)!.

\item[$\circ$] \lstinline!writeKiMoKi(outfile)!

Writes the definition of the main algorithm of the solver module,
\lstinline!kiMoKi(z)!. The routine \lstinline!kiMoKi(z)! processes
the \emph{kick}-\emph{push}-\emph{move}-\emph{kick} substeps of
a full timestep, including the iterative solution of equation
\eqref{CanonicalTransformation}. 

\item[$\circ$] \lstinline!writeRunExampleDP(outfile, modname, qvars, pvars, params)!,\\
\lstinline!writeRunExampleMP(outfile, modname, qvars, pvars, params)!
 
Writes a simple example program illustrating how to use the solver module.
Some routines of this program solves the Hamilton's equation over a time interval,
with random initial conditions and parameters (which most likely must first be manually
changed to sensible values), and writes a plot of the solution
to a \lstinline!pdf! file. Other routines check how well energy is conserved
by the solver, for a set of timesteps, and writes a plot of the energy errors
to another \lstinline!pdf! file.

\end{itemize}

\begin{figure}[H]
\label{exampleUseOfCodeGenerator}
\begin{center}
\includegraphics[width=0.99\textwidth]{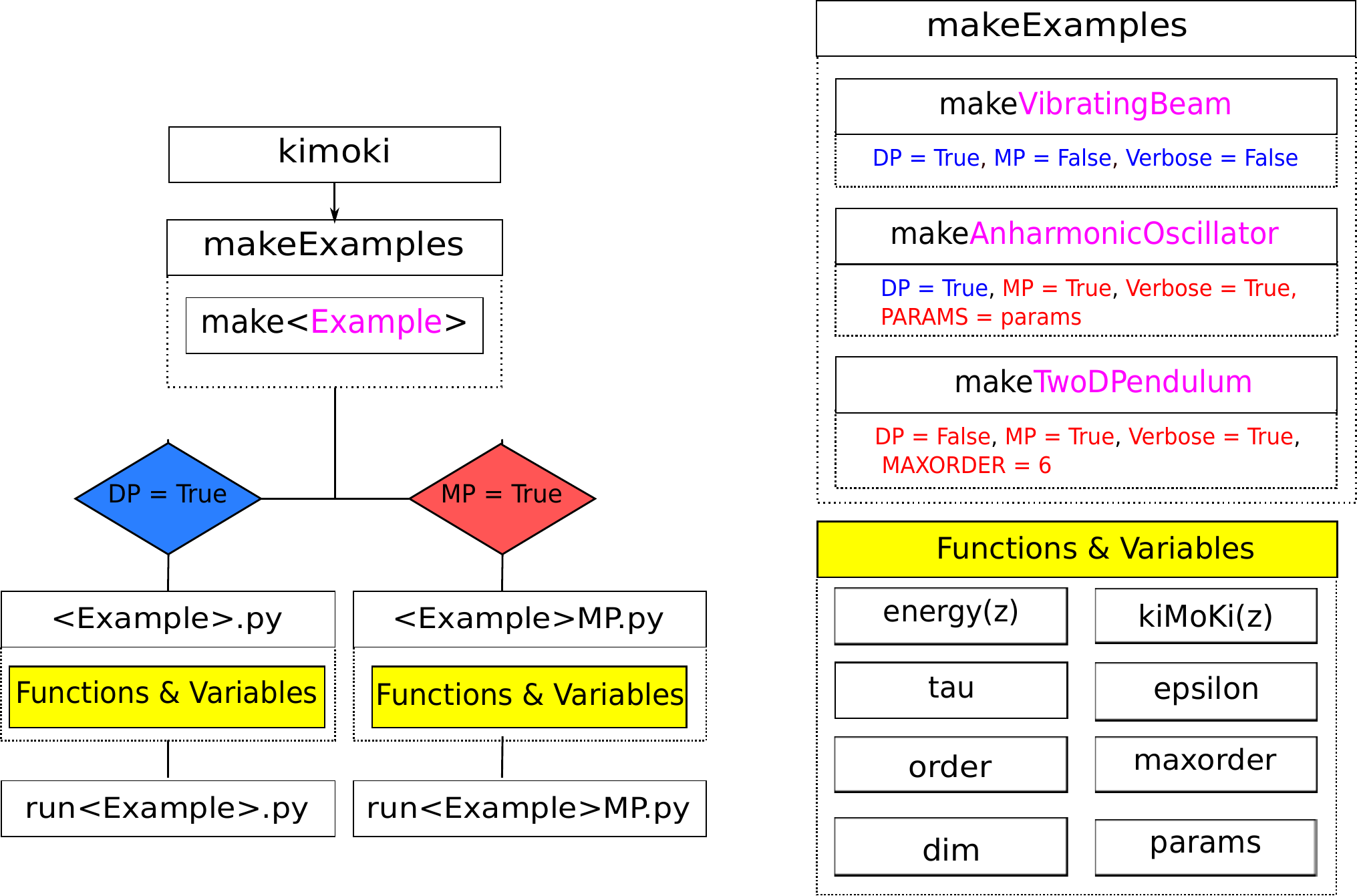}
\end{center}
\caption{This figure illustrates use of the code generator. The call of 
\lstinline!kimoki.makeModules('<Example>',...)! will generate a solver
module, \texttt{<Example>.py}, and a demonstration runfile 
\texttt{run<Example>.py}. When the optional parameter \lstinline!MP=True!
a multiprecision version of the solver module and runfile is generated.
The solver module consists of various functions and variables. Its most
important function is \lstinline!kiMoKi(z)!, which updates the solution
\lstinline!z! through one full timestep. The function \lstinline!energy(z)!
evaluates the Hamiltonian at the phase space point \lstinline!z!. 
Many other functions are also defined. F.i., \lstinline!T2(z)!, \lstinline!T4(z)!, 
\lstinline!T6(z)! (these return 0 if \lstinline!ALLTK=False!), and
\lstinline!V2(z)!, \lstinline!V4(z)!, \lstinline!V6(z)!. Plus additional
functions which are not intended to be called directly by the user. Several
parameters, some which can be changed by the user, are also defined: The timestep
\lstinline!tau!, the accuracy \lstinline!epsilon! to which equation
\eqref{CanonicalTransformation} must be solved, which \lstinline!order!
of the integrator to use when running \lstinline!kiMoKi(z)!. The maximum
order \lstinline!maxorder! of solvers available (must not be changed by the user),
the phase space dimension \lstinline!dim! of the model being solved 
(must not be changed by the user),
and a (possibly empty) list of parameters \lstinline!params! on which
the Hamiltonian depends (must not be changed by the user).}
\end{figure}

\begin{figure}[H]
\label{runexampleUseOfCodeGenerator}
\begin{center}
\includegraphics[width=0.99\textwidth]{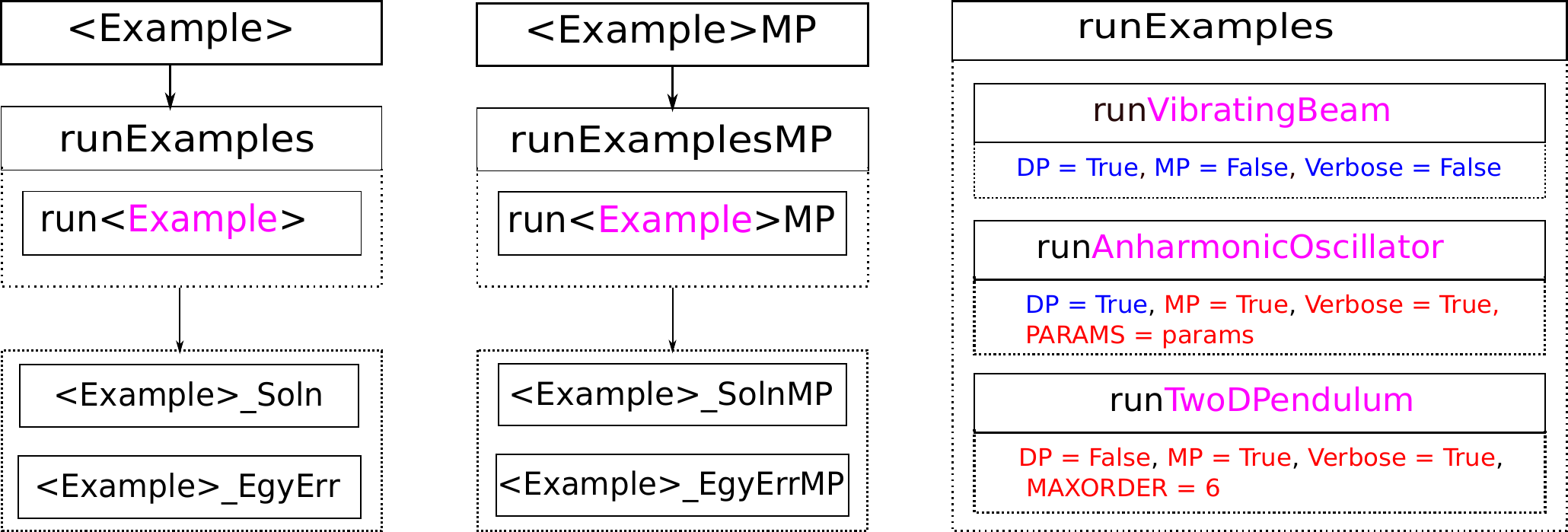}
\end{center}
\caption{This figure illustrates the runfile example. Normally the runfile
generates two figure files in \texttt{.png} format (this can easily be changed
to \texttt{.pdf} format by the user).
The runfile also generate several intermediate
files. Their names begin with the symbol \texttt{\#}. These files are normally
deleted after use (this can easily be changed by the user).}
\end{figure}

\section{Concluding remarks}\label{concludingRemarks}

In this paper we have demonstrated that the proposed extensions
of the standard St{\"o}rmer-Verlet symplectic integration scheme
can be implemented numerically, and that the implemented code behave
as expected with respect to accuracy. Here we have not focused on time or
memory efficiency of the generated code, which may be viewed as a
reference implementation known to work correctly. We have experienced
this to be a good starting point for manual implementation of
more efficient code for large, structured systems, f.i. the
Fermi-Pasta-Ulam-Tsingou type lattice models studied in \cite{AAK11},
and molecular dynamics type simulations studied in \cite{AAKR}.
Code for the latter systems are quite straightforward to implement
using NumPy arrays \cite{WCV}, which also leads to efficient working code.

It is also straightforward to modify our program to
generate code in other computer languages.

\section*{Acknowledgements}

We thank professor Anne Kv{\ae}rn{\o} for useful discussions,
helpful feedbacks, and careful proofreading. We also acknowledge support
provided by Statoil via Roger Sollie, through a professor II grant
in Applied mathematical physics.






\end{document}